\documentclass[10pt,amsmath,amssymb,amsfonts]{article}
\usepackage{epstopdf}
\usepackage{epsfig}
\usepackage{graphicx} 

\newcommand{\be}{\begin{equation}}
\newcommand{\ee}{\end{equation}}
\newcommand{\bea}{\begin{eqnarray}}
\newcommand{\eea}{\end{eqnarray}}

\begin{document}

\title{Stochastic mRNA production by a three-state gene}
\author{Boris Rubinstein, Jay Unruh, Julia Zeitlinger
\\Stowers Institute for Medical Research
\\1000 50$^{}\mbox{th}$ St., Kansas City, MO 64110, U.S.A.}
\date{\today}
\maketitle

\begin{abstract}
We consider a model of mRNA production 
governed by the dynamics of a gene that exists in three possible states --
inactive, poised and active. The transitions between the adjacent states
are controlled by stochastic processes characterized by 
corresponding on/off rates.
mRNA is produced only when the gene is in active state and
we also consider mRNA denaturation leading to its death.
We derive the distribution of mRNA number and compare it
to the known result for the two-state gene model.
\end{abstract}

\section{Introduction}

Consider a model of mRNA synthesis by a gene having three states -- inactive (0), 
poised/paused/waiting (1) and active (2).
mRNA translation takes place when the gene is in active state only and it degrades 
with a rate proportional to its 
concentration (number). Our approach extends the one used in \cite{Pessoud1995} for a two-state gene
that allows to compute the probability distribution of mRNA number as a function of the rates. 
An alternative way to describe gene state occupancy dynamics is based on explicit 
inclusion of the noise term into dynamical equations \cite{Rieckh2014} which
focused on the inflence of noise on the mRNA production.

\begin{figure*}[h]
\begin{center}
\psfig{figure=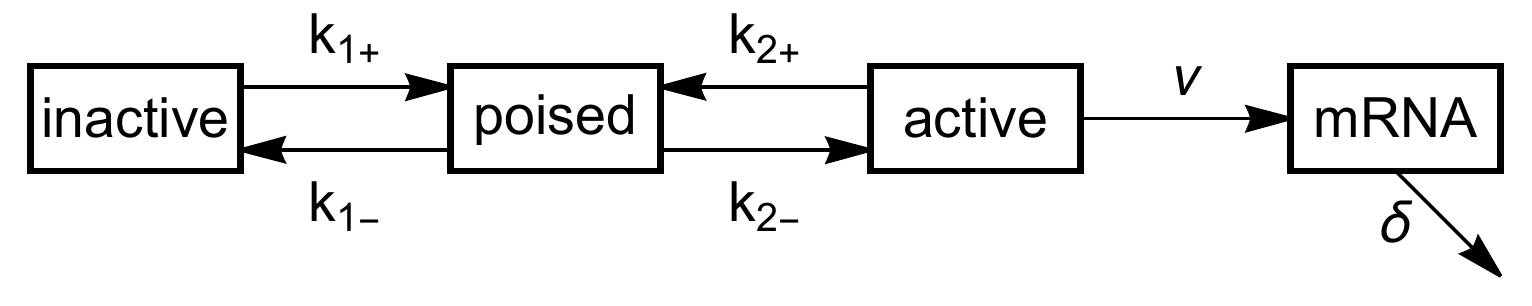,width=12cm}
\end{center} 
\caption{The transitions in the three-state gene model.}
\label{Scheme}
\end{figure*}

The state of the system (gene + mRNA) is described by a two-dimensional integer
vector $\{g,n\}$ where $g=0,1,2$ is the gene state and $n$ denotes a number of mRNA molecules.
The transitions between gene states (shown in Fig.\ref{Scheme}) occur with the rates $k_{1\pm}$ for the pair inactive--poised and 
$k_{2\pm}$ for the pair poised--active while the direct transitions between inactive and active states are
forbidden. The gene in active state produces a mRNA molecule with the rate $\nu$, these molecules
degrade with the rate $\delta$.

All allowed transitions are described below.
\vskip0.3cm
$
\begin{array}{ll}
\mbox{Transition} & \mbox{Rate} \\
\{0,n\} \to \{1,n\} &  k_{1+} \\
\{0,n\} \to \{0,n-1\} &  \delta n \\
 & \\
\{1,n\} \to \{0,n\} &  k_{1-} \\
\{1,n\} \to \{1,n-1\} &  \delta n \\
\{1,n\} \to \{2,n\} &  k_{2+} \\
 & \\
\{2,n\} \to \{1,n\} &  k_{2-} \\
\{2,n\} \to \{2,n-1\} &  \delta n \\
\{2,n\} \to \{2,n+1\} &  \nu 
\end{array}
$
\vskip0.3cm

\section{Balance equations}
Using the above characterization of the transitions we can write down the 
dynamical equations  for the probability $p_{g,n}(t)$ to have $n$ copies of mRNA
for a gene in state $g$ at time $t$
\bea
p_{0,n}'(t) &=& -k_{1+} p_{0,n} -n \delta p_{0,n} + (n+1)\delta 
p_{0,n+1} + k_{1-} p_{1,n}, 
\nonumber 
\\
p_{1,n}'(t) &=& -k_{1-} p_{1,n} -n \delta p_{1,n}+ (n+1)\delta 
p_{1,n+1} + k_{1+} p_{0,n} -k_{2+} p_{1,n} + k_{2-} p_{2,n},
\nonumber\\
p_{2,n}'(t) &=& -k_{2-} p_{2,n}+ k_{2+} p_{1,n}-n \delta p_{2,n} + (n+1)\delta 
p_{2,n+1}-\nu p_{2,n}+\nu p_{2,n-1}.
\nonumber
\eea
Rescaling the time variable $t \to t \delta$ and all rates $r \to r/\delta$ we obtain
\bea
p_{0,n}'(t) &=& -k_{1+} p_{0,n} -n p_{0,n} + (n+1)\delta 
p_{0,n+1} + k_{1-} p_{1,n}, 
\label{eqP}
\\
p_{1,n}'(t) &=& -k_{1-} p_{1,n} -n p_{1,n}+ (n+1)
p_{1,n+1} + k_{1+} p_{0,n} -k_{2+} p_{1,n} + k_{2-} p_{2,n},
\nonumber\\
p_{2,n}'(t) &=& -k_{2-} p_{2,n}+ k_{2+} p_{1,n}-n  p_{2,n} + (n+1) 
p_{2,n+1}-\nu p_{2,n}+\nu p_{2,n-1}.
\nonumber
\eea

\section{Generating functions}
The probabilities $p_{i,n}(t)$ give rise to generating functions $G_i(z,t)$ defined
as follows
\be
G_i(z,t) = \sum_{n=0}^{\infty}
z^n p_{i,n}(t).
\label{Gi_def}
\ee
Setting above $z=1$ we find 
$$
G_i(1,t) = \gamma_i = \sum_{n=0}^{\infty} p_{i,n}(t)
$$
the probability of a gene to be in $i$-th state.
The only observable of the model is the mRNA number $n$ at time $t$ 
defined by a probability $p_n =  p_{0,n}+ p_{1,n}+ p_{2,n}$. The corresponding
generating function reads 
$$
G(z,t) = G_0(z,t) + G_1(z,t)  + G_2(z,t)  .
$$
Note that
$$
\partial G_i(z,t)/\partial z = 
 \sum_{n=0}^{\infty} n z^{n-1} p_{i,n}(t) = 
 \sum_{n=0}^{\infty} (n+1) z^{n} p_{i,n+1}(t),
$$
leading to 
$$
 \sum_{n=0}^{\infty} n z^{n} p_{i,n}(t)  = z\; \partial G_i(z,t)/\partial z,
$$
while (\ref{Gi_def}) implies 
$$
 \sum_{n=0}^{\infty} z^{n} p_{i,n-1}(t) = z G_i(z,t).
$$
In order to derive equations for the generating functions we 
multiply each equation in (\ref{eqP}) by $z^n$ and sum up
w.r.t. $n$. This procedure leads to 
\bea
\frac{\partial G_0(z,t)}{\partial t} &=& 
-k_{1+} G_0(z,t)+k_{1-} G_1(z,t)+(1-z) \frac{\partial G_0(z,t)}{\partial z},
\label{eqG}
\\
\frac{\partial G_1(z,t)}{\partial t} &=& 
-k_{1-} G_1(z,t)+k_{1+} G_0(z,t)-k_{2+} G_1(z,t)+k_{2-} G_2(z,t)
+(1-z) \frac{\partial G_1(z,t)}{\partial z},
\nonumber \\
\frac{\partial G_2(z,t)}{\partial t} &=& 
-k_{2-} G_2(z,t)+k_{2+} G_1(z,t)+(1-z) \frac{\partial G_2(z,t)}{\partial z}
-\nu(1-z)\frac{\partial G_2(z,t)}{\partial z}.
\nonumber
\eea
The natural assumption in the model is that at $t=0$ the 
gene is inactive ($i=0$) and the number of mRNA molecules $n=0$, 
so that the initial conditions read $p_{0,0}(0)=1$ and we find
$G_i(z,0)=\delta_{i,0}$.
Adding up the equations in (\ref{eqG}) we obtain
\be
\frac{\partial G(z,t)}{\partial t} = 
(1-z) \left[\frac{\partial G(z,t)}{\partial z}
-\nu\frac{\partial G(z,t)}{\partial z}\right].
\label{const0}
\ee

\section{Steady state equations}
Analytical solution of (\ref{eqG}) that completely determines
the dynamics of the mRNA is not known. We consider an asymptotic
behavior of the probability $p_{i,n}$ at large times $t \to \infty$.
Then the equations (\ref{eqG}) turn into a system of ODEs
\bea
&&
k_{1-} G_1(z)-k_{1+} G_0(z)+(1-z) G_0'(z)
=0,
\label{eqG0}
\\
&& 
k_{1+} G_0(z)-k_{1-} G_1(z)-k_{2+} G_1(z)+k_{2-} G_2(z)
+(1-z) G_1'(z)=0,
\label{eqG1}\\
&& 
k_{2+} G_1(z)-k_{2-} G_2(z)+(1-z)G_2'(z)
-\nu(1-z)G_2'(z)=0.
\label{eqG2}
\eea
Adding up these equations we obtain a relation
\be
\sum_{i=0}^2 G_i'(z) = G'(z) = \nu G_2(z).
\label{const1}
\ee
In order to find the asymptotic values of $\gamma_i$
set $z=1$ in (\ref{eqG0}-\ref{eqG2}) and obtain
\bea
&&
k_{1-} \gamma_1-k_{1+} \gamma_0=0,
\nonumber
\\
&& 
k_{1+}\gamma_0-k_{1-}\gamma_1-k_{2+} \gamma_1+k_{2-}\gamma_2=0,
\nonumber\\
&& 
k_{2+} \gamma_1-k_{2-} \gamma_2=0,
\nonumber
\eea
where $\gamma_0+\gamma_1+\gamma_2=1$
to produce 
\be
\gamma_0=k_{1-}k_{2-}/K,
\quad
\gamma_2=k_{1+}k_{2+}/K,
\quad
K =k_{1-}k_{2-}+k_{1+}k_{2-}+k_{1+}k_{2+}.
\label{gamma}
\ee

\section{Solution for $G_2(z)$}
It follows from (\ref{const1}) that 
in order to find $G(z)$ and determine the observable probability $p_n$ 
 it is suffice to have an explicit expression for $G_2(z)$.  
In this section we first reduce the system (\ref{eqG0}-\ref{eqG2}) to
a single equation for $G_2(z)$ and then obtain its solution.

\subsection{Reduction to a single equation} 

First eliminate $G_1$ from (\ref{eqG0},\ref{eqG2}) to produce
\be
k_{1-}[(k_{2-}+\nu(1-z))G_2-(1-z) G_2'] = 
k_{2+}[k_{1+}G_0-(1-z) G_0'].
\label{G0G2a}
\ee
Differentiate (\ref{eqG0}) w.r.t. $z$ and substitute into the 
result $G_1'(z)$  from (\ref{const1})
$$
(1-z)G_0''-(1+k_{1+})G_0'+k_{1-}(\nu G_2-G_0'-G_2')=0.
$$
Now use this relation together with (\ref{G0G2a}) to obtain
\be
(1-z)(G_2''-\nu G_2')-(1+k_{2-}+k_{2+})G_2'+\nu(1+k_{2+})G_2=k_{2+}G_0'.
\label{G0G2b}
\ee
Similarly differentiate (\ref{eqG2}) w.r.t. $z$ and substitute $G_1'(z)$ 
to produce
$$
(1-z)[G_2'-\nu G_2']'-(1+k_{2-}G_2')+\nu G_2+k_{2+}(\nu G_2-G_0'-G_2')=0.
$$
Use it again with (\ref{G0G2a}) to generate
\be
-(1-z)G_0''-(1+k_{1-}+k_{1+})G_0'=k_{1-}(\nu G_2-g_2').
\label{G0G2c}
\ee
Finally, express $G_0'$ from (\ref{G0G2b}) and use it in (\ref{G0G2b})
to construct a single third order differential equation for $G_2$
\bea
&& (z-1)^2 G_2''' + (z-1)(3+\kappa_1-\nu(z-1))G_2''-
\nu  \kappa_3    G_2
\nonumber \\
&+&
[\kappa_3+ (1+k_{1-}+k_{1+})k_{2-}-\nu(z-1)(3+\kappa_2)]
G_2'=0,
\label{eqG2only}
\eea
where
$$
\kappa_1 = k_{1-}+k_{1+}+k_{2-}+k_{2+},
\kappa_2 = \kappa_1 - k_{2-},
\kappa_3 = k_{1-}+(1+k_{1+})(1+k_{2+}).
$$

\subsection{Expression for $G_2(z)$}
The solution of the linear ODE (\ref{eqG2only}) was obtained
with computer algebra software {\it Mathematica},
it has three components, but only one of these three
does not diverge at $z=1$. This component
is a generalized hypergeometric function
$$
G_2(z) = c\;
{}_2F_2(\{1+K_{2-},1+K_{2+}\},\{1+K_{1-},1+K_{1+}\},\nu(z-1)),
$$
where $c$ is the undetermined constant and we introduce the shortcut notations
$$
K_{1\pm} = (\kappa_1 \pm \sqrt{\kappa_1^2-4\kappa_0})/2,
\quad
K_{2\pm} = (\kappa_2 \pm \sqrt{\kappa_2^2-4k_{1+}k_{2+}})/2,
$$
where $\kappa_0 = (k_{1-}+k_{1+})k_{2-}+k_{1+}k_{2+}$.
Note that $K_{i\pm}$ are homogeneous functions of degree one. 

The constant $c$ is selected to satisfy the condition
$G_2(1) = \gamma_2$. 
As ${}_2F_2(\{a_1,a_2\},\{b_1,b_2\},0)=1$ we immediately find that 
$c=\gamma_2= k_{1+}k_{2+}/(k_{1-}k_{2-}+k_{1+}k_{2-}+k_{1+}k_{2+})$.
Thus we finally find
\be
G_2(z) =
\frac{
{}_2F_2(\{1+K_{2-},1+K_{2+}\},\{1+K_{1-},1+K_{1+}\},\nu(z-1))}
{1 + (1+k_{1-}/k_{1+})(k_{2-}/k_{2+})}.
\label{G2sol}
\ee

\section{Computation of $G(z)$}
Relation (\ref{const1}) allows to find $G(z)$ as
an integral
$$
G(z) = \nu \int G_2(z) dz, 
$$
provided $G(1) = 1$.
The  generalized hypergeometric functions ${}_pF_q({\bf a},{\bf b},x)$, where
${\bf a}= \{a_1,a_2,\ldots,a_p\}$ and ${\bf b}= \{b_1,b_2,\ldots,b_q\}$ are the vectors of $p$ and $q$ components respectively,
have a nice property that both derivatives and integrals are expressed through
the same functions 
\cite{WolframFunc1}. Specifically,
$$
\int {}_pF_q({\bf a},{\bf b},x) dx = 
\frac{\pi({\bf b}-1)}{\pi({\bf a}-1)}\;
{}_pF_q({\bf a}-1,{\bf b}-1,x),
$$
where we employ the following notations for a $p$-dimensional vector ${\bf v}$
$$
\pi({\bf v}) = \prod_{i=1}^p v_i,
\quad
{\bf v}+m = \{v_1+m,v_2+m,\ldots,v_p+m\}.
$$
Introducing $x=\nu(z-1)$ we find 
$$
\nu \int {}_pF_q({\bf a},{\bf b},\nu(z-1)) dz= 
\int {}_pF_q({\bf a},{\bf b},x) dx,
$$
so that
\be
G(z) =
{}_2F_2(\{K_{2-},K_{2+}\},\{K_{1-},K_{1+}\},\nu(z-1)),
\quad
G(1) = 1.
\label{Gsol}
\ee

\section{Steady state mRNA distribution}
From the definition of the generating function we have
$$
G(z) = \sum_{n=0}^{\infty}
z^n p_{n},
$$
leading to a conclusion that $p_n$ is the 
coefficient in the Taylor expansion of $G(z)$ around $z=0$.
The explicit expression for $p_n$ then reads \cite{WolframFunc2}
\bea
p_n &=& 
\frac{r_n(K_{2-})r_n(K_{2+})}{r_n(K_{1-})r_n(K_{1+})}\cdot 
\frac{\nu^n}{n!} \cdot \; {}_2F_2(\{K_{2-},K_{2+}\}+n,\{K_{1-},K_{1+}\}+n,-\nu),
\nonumber \\
r_n(x) &=&  \Gamma(x+n)/\Gamma(x),
\label{pn}
\eea
where $r_n(x)$ is the Pochhammer symbol defined via
gamma function $\Gamma(x)$.

Note that for $x \sim O(1)$ and large $n\gg 1$ 
the value of $r_n(x)$ grows as $r_n(x) \sim x^n$
so that the first factor in (\ref{pn}) behaves as 
$(K_{2-}K_{2+})^n/(K_{1-}K_{1+})^n$.

\subsection{Reduction and comparison to two-state gene model}
To reduce the three-state model to two-state one we set $k_{1-}=0$ so that $\gamma_0 = 0$ and
the gene can be only in the active (2) or the poised (1)  state which effectively plays a role of the 
inactive state. In this case $p_{0,n}=0,\; G_0(z)=0$ and we end up with the following system of equations
\bea
&&
-k_{2+} G_1(z)+k_{2-} G_2(z)
+(1-z) G_1'(z)=0,
\label{eqG1twost}\\
&& 
k_{2+} G_1(z)-k_{2-} G_2(z)+(1-z)G_2'(z)
-\nu(1-z)G_2'(z)=0.
\label{eqG2twost}
\eea
The model (\ref{eqG1twost},\ref{eqG2twost}) was discussed in \cite{Pessoud1995}
and in this case the probability $p_n$ reads
\be
p_n = 
\frac{\nu^n r_n(K_{+})}{n! r_n(K_{+}+K_{-})}\;
{}_1F_1(\{K_{+}+n\},\{K_{+}+K_{-}+n\},-\nu),
\label{pn2st}
\ee
where $K_{\pm} = k_{2\pm}$.
This result can be also obtained directly from (\ref{pn}).
We find for $k_{1-} = 0$ the values of 
$\kappa_0 = k_{1+}(k_{2-}+k_{2+}) ,\ \kappa_1 = k_{1+}+k_{2+}+k_{2-}, \ 
\kappa_2 =  k_{1+}+k_{2+}$ to compute
$K_{1\pm}=( k_{1+}+k_{2+}+k_{2-} \pm( k_{1+}-k_{2+}-k_{2-}))/2$ 
and $K_{2\pm} = ( k_{1+}+k_{2+} \pm( k_{1+}-k_{2+}))/2$.
It leads to $K_{1-}=K_{2+} = k_{1+}$ and $K_{1+}=k_{2+}+k_{2-}=K_{+}+K_{-}, \ 
K_{2-} =k_{2+}=K_{+}$ so that the property of 
the hypergeometric function 
$$
{}_{p+1}F_{q+1}(\{a_1,a_2,\ldots,a_p,c\}, \{b_1,b_2,\ldots,b_q,c\},x) = 
{}_{p}F_{q}(\{a_1,a_2,\ldots,a_p\}, \{b_1,b_2,\ldots,b_q\},x)
$$
implies
$$
{}_2F_2(\{K_{2-},K_{2+}\}+n,\{K_{1-},K_{1+}\}+n,-\nu) = 
{}_1F_1(\{K_{+}+n\},\{K_{+}+K_{-}+n\},-\nu).
$$

\subsection{Qualitative model predictions}
\begin{figure}[h]
\begin{tabular}{cc}
\includegraphics[width=7.5cm]{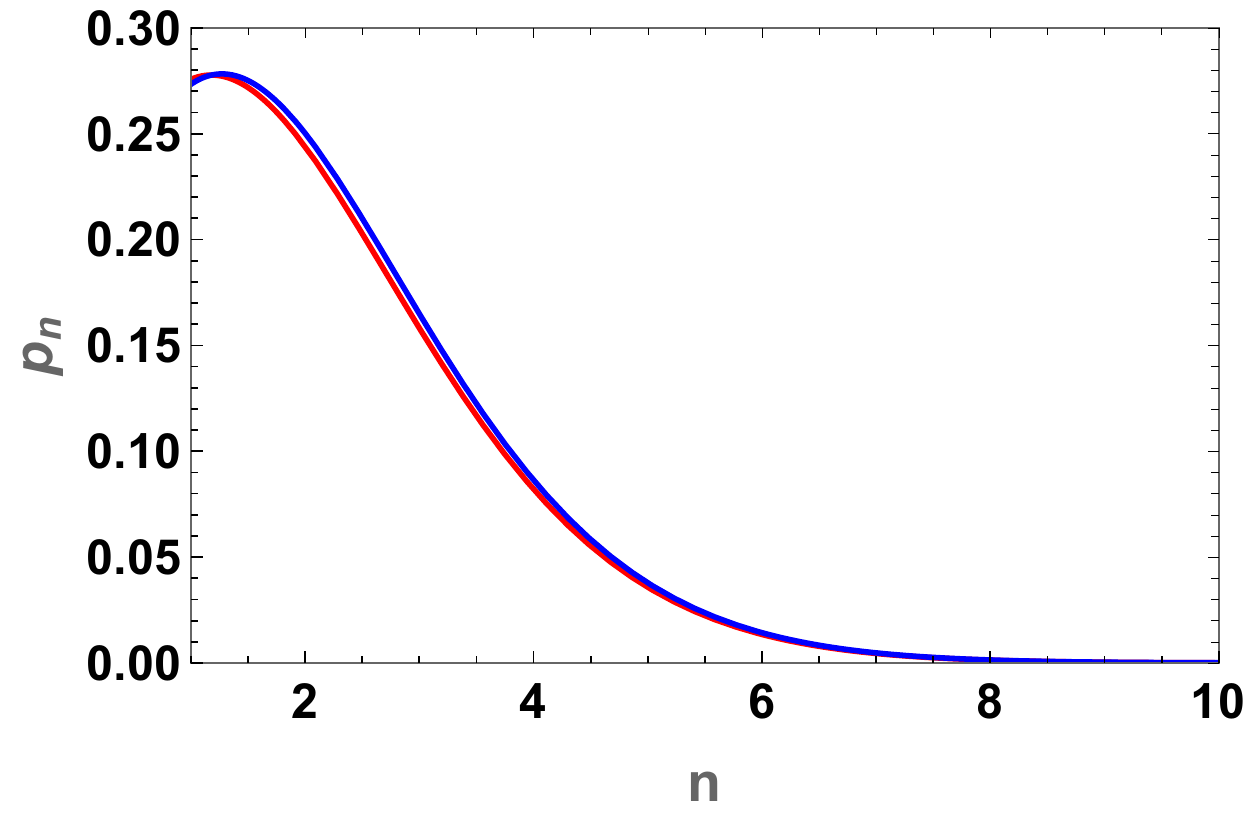} &
\includegraphics[width=7.5cm]{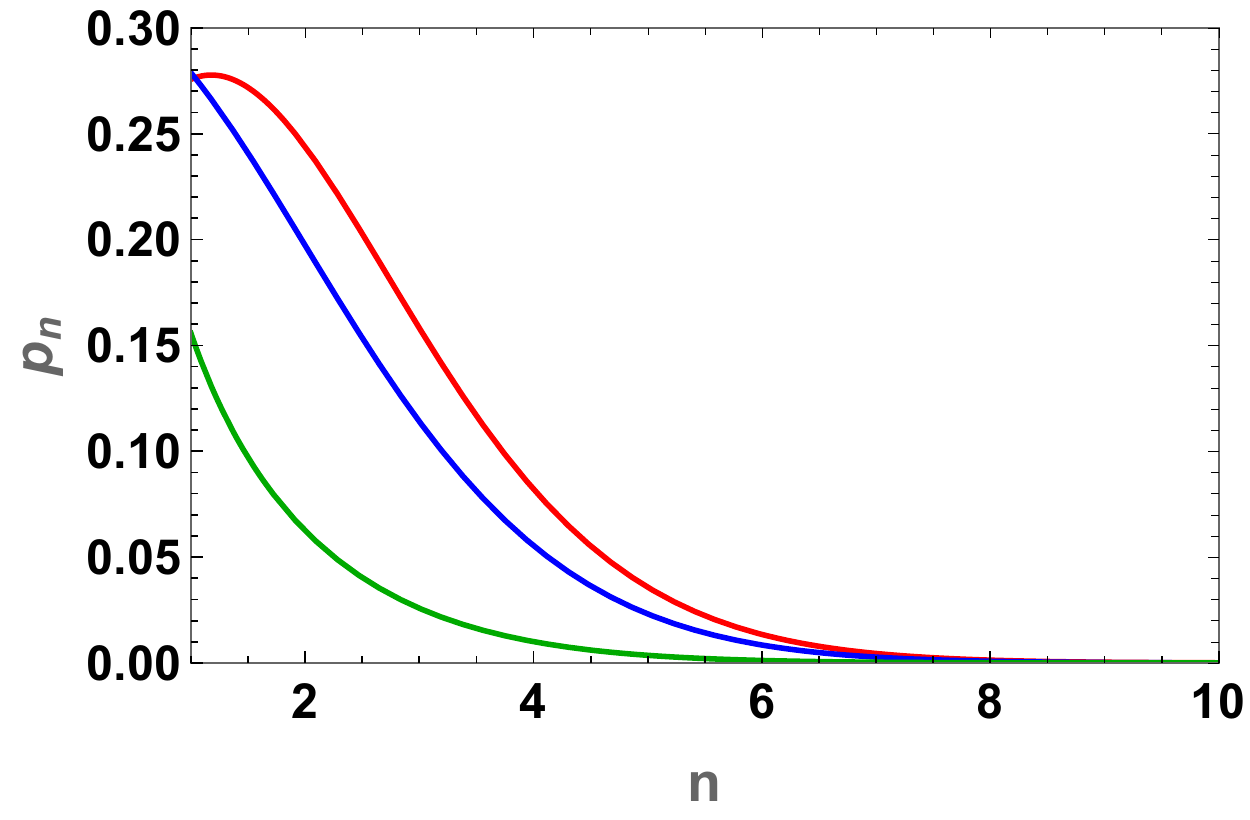} \\
(a) & (b) 
\end{tabular}
\caption{(a) The two-state (blue) and three-state (red) gene
model mRNA distributions for $k_{1-} \ll k_{1+}$ demonstrate 
nearly perfect coincidence. The parameter values are
$k_{1-}=0.13,\; k_{1+}=1.3,\; k_{2-}=2.3,\; k_{2+}=4.2,\; \nu=3$ for
the three-state and $k_{-}=2.3,\; k_{+}=4.2,\; \nu=3$ for the two-state gene
model.
(b) An effect of the inactive state population on mRNA production in three-state
gene model. The parameters are
$k_{1+}=1.3,\; k_{2-}=2.3,\; k_{2+}=4.2,\; \nu=3$
with $k_{1-}$ equal to $0.13$ (red), $1.3$ (blue) and $13.0$ (green)
respectively.}
\label{fig1}
\end{figure}

When $k_{1-} \ll k_{1+}$ one expects that $p_n$ for three-state gene
should be quite close to the two-state gene solution. The numerical simulations
confirm the assumption (Fig. \ref{fig1}a). 
When $k_{1-}$ increases the
population of the inactive state according 
to (\ref{gamma}) also grows
thus reducing simultaneously the active state probability.
As expected this reduction results in the
shift to the cells expressing lower mRNA numbers (Fig. \ref{fig1}b).

\section{Hypergeometric function asymptotic expansion for large argument}
The $p_n$ expressions 
(\ref{pn2st}) and (\ref{pn}) for two- and three-state model respectively
might need to be evaluated when both the argument $\nu$ and the parameters of hypergeometric function
are large. In this case it is worth to use asymptotic expansion for fast and
accurate computation of these functions ${}_kF_k,\; k=1,2$.

\subsection{Asymptotics of ${}_1F_1$}
We use  asymptotics at large $|z| \to \infty$ 
\be
{}_1F_1(\{a\},\{b\},z) \sim 
\left\{
\begin{array}{ccc}
(-z)^{-a} \Gamma(b)
{}_2F_0(\{a,a-b+1\},\{\},-1/z)/\Gamma(b-a), 
&\quad & a,b < |z|,
\\
z^{a-b} e^z \Gamma(b)
{}_2F_0(\{b-a,1-a\},\{\},1/z)/\Gamma(a), 
&\quad & a,b \ge |z|,
\end{array}
\right.
\label{1F1_asymp}
\ee
where the function ${}_2F_0(\{a,b\},\{\},x)$ is computed through 
the series 
$$
{}_2F_0(\{a,b\},\{\},x) = 
\sum_{k=0}^{\infty} r_k(a) r_k(b) x^k/k!,
$$
as a particular case of the general definition 
\be
{}_{p}F_{q}(\{a_1,a_2,\ldots,a_p\}, \{b_1,b_2,\ldots,b_q\},x) =
\sum_{k=0}^{\infty} 
\frac{\prod_{i=1}^{p} r_k(a_i)}{\prod_{j=1}^{q} r_k(b_j)}\cdot  \frac{x^k}{k!},
\label{pFq_general}
\ee
where the upper limit in the sum is replaced by positive integer $m$ representing number of 
terms in the expansion.
Note that (\ref{pn2st}) for $K_{\pm} \sim O(1)$ at large $\nu \gg 1$ and $n < \nu$ one has to use 
the first formula in (\ref{1F1_asymp}) while for large $n > \nu$ the second line is employed.

\subsection{Asymptotics of ${}_2F_2$}
When $|z| \gg 1$ we employ the following asymptotic formula \cite{WolframFunc3}
\bea
{}_2F_2(\{a_1,a_2\},\{b_1,b_2\},z) \sim 
\left\{
\begin{array}{lc}
(-z)^{-a_1} \frac{\Gamma(b_1)\Gamma(b_2)\Gamma(a_2-a_1)}
{\Gamma(a_2)\Gamma(b_1-a_1)\Gamma(b_2-a_1)}&\\
\times
{}_3F_1(\{a_1,a_1-b_1+1,a_1-b_2+1\},\{a_1-a_2+1\},-1/z)
& \\
+ 
(-z)^{-a_2} \frac{\Gamma(b_1)\Gamma(b_2)\Gamma(a_1-a_2)}
{\Gamma(a_1)\Gamma(b_1-a_2)\Gamma(b_2-a_2)}&\\
\times
{}_3F_1(\{a_2,a_2-b_1+1,a_2-b_2+1\},\{a_2-a_1+1\},-1/z), &
a_i,b_i < |z|,
\\
e^z z^{A_2-B_2} \frac{\Gamma(b_1)\Gamma(b_2)}
{\Gamma(a_1)\Gamma(a_2)}
\sum_{k=0}^{\infty} c_k (1/z)^{k},
& a_i,b_i \ge|z|,
\end{array}
\right.
\nonumber 
\eea
where the upper limit in the last line sum is replaced by positive integer $m$ .
We use the notation $A_2=a_1+a_2,\; B_2=b_1+b_2$ and the 
expansion coefficients $c_k$ are defined by a recursion
\bea
&&c_0=1, \
c_1 = (A_2-1)(A_2-B_2)+b_1b_2-a_1a_2,\
\nonumber\\
&&
k c_k = (1 - B_2 + a_1(2 +a_1) + a_2(2 +a_2) - A_2 B_2 + a_1 a_2 + b_1 b_2 + 
(2 B_2 - 3 (A_2 + 1)) k + 2 k^2) c_{k-1} 
\nonumber \\
&&- (k - A_2 + b_1 - 1) (k - A_2 + b_2 - 1) (k - A_2+B_2 - 2) c_{k-2}.
\label{coeff_ck}
\eea

\end{document}